\documentclass[article,oneside,12pt]{article}

\input{epsf}
\usepackage{amsfonts,amsmath,amssymb,amsthm,mathrsfs}
\usepackage[dvips]{graphics}

\theoremstyle{plain}
\newtheorem{theo}{Theorem}
\newtheorem{lemo}{Lemma}

\theoremstyle{remark}
\newtheorem{rema}{Remark}

\begin{document}
\title{On spatial thinning-replacement processes based on Voronoi cells }

\author{K.A. Borovkov\thanks{Department of Mathematics and Statistics, The University of Melbourne, Parkville 3010, Australia
}, \ \ D.A. Odell\thanks{MASCOS, 139 Barry Street, Carlton, Melbourne, Australia 3010}}

\maketitle

\def\vex#1{\textbf{\textit{#1}}}
\def\tdots{,\dots ,}
\def\muN{\mu^{\otimes N}}
\def\Pr{\text{Pr}}

\begin{abstract}
We introduce a new class of spatial-temporal point processes based on Voronoi
tessellations.  At each step of such a process, a point is chosen at random according
to a distribution determined by the associated Voronoi cells. The point is then
removed, and a new random point is added to the configuration. The dynamics are simple
and intuitive and could be applied to modeling natural phenomena. We prove ergodicity
of these processes under wide conditions.
\\



\end{abstract}

\section{Introduction}

A spatial point process is a stochastic model for the location of events in a space,
and as such it is a random element which takes discrete sets of points as its values.
Trees in a forest, schools in a city, or capillaries on the surface of a bodily organ
are examples of situations that can be modelled by a point process over a
two-dimensional manifold, impurities in metals or positions of submarines require
models defined over three-dimensional manifolds. The one-dimensional point process is
often used as a model for events in time, but may also be applied to such problems as
the physical distribution of files on a hard-disk where the underlying space can be
treated as an interval in $\mathbb{R}$.

The Voronoi tessellation is a useful tool for the analysis of spatial point processes.
This can be defined for any locally finite set of points $X$ in a metric space $(M,d)$.
With each $x\in X$ we associate its \emph{Voronoi cell} \begin{equation}C_{x}^{X}(M):=
\left\{y \in M : d(y,x)=\min{[d(y,x'):x'\in X]}\right\}.\end{equation} [We will
generally suppress the reference to $M$ in the notation for the cell and just write
$C_{x}^{X}$ where the underlying set is clear from the context.] The
\textit{tessellation} produced by $X$ is $\mathscr{T}_{X}=\left\{C_{x}^{X}:x \in
X\right\}$, and the elements of $X$ are called the \textit{generators} of
$\mathscr{T}_{X}$. Statistics drawn from analysis of this structure provide an
intrinsic description of the distribution pattern of the set of generators and have
accordingly been investigated extensively (see e.g. \cite{OBS}). The Voronoi
tessellation is also a natural object of interest whenever concepts such as
\textit{catchment area} or \textit{zone of influence} are appropriate to the situation
being modelled, as in most of the examples mentioned above. When this is the case it is
also very interesting to consider point processes \textit{evolving in time} under laws
which are determined by characteristics of the Voronoi tessellation.  A simple example
is a process in which the lifetime of a point is a random variable with distribution
determined by its Voronoi cell. Examples of this kind of spatial-temporal point process
are the Hotelling processes \cite{OS87}, which model the geographical distribution of
businesses competeing for market share and the Adjustment models for territorial animal
behaviour \cite{HT76}. As far as we know there has been no general treatment of such
models although they promise to produce a rich class of stochastic objects typifying
various kinds of spatial clustering.

In this paper we  consider some classes of such models viewed as discrete-time Markov
processes $\left\{\mathbf{X}_n\right\}_{n\geq 0}$ taking values in a fixed
finite-dimensional configuration space. Let $(M,d)$ be a metric space and $N$ a fixed
positive integer representing the initial number of points in $M$. The associated
configuration space will be either  $\mathscr{X}=\bigcup_{k=1}^{N}M^{k}$ or
$\mathscr{X}=M^{N}$, depending on whether the process is of a thinning or
thinning-replacement type. Any finite set of points $X$ generating a Voronoi
tessellation can of course be represented in many different ways as a point $\vex{x}\in
\mathscr{X}$ since the order of the points is immaterial to the tessellation, however
it is convenient to retain an ordering of the points in the configuration space. In the
thinning-replacement process the total number of points after each replacement remains
constant, although the thinning and replacement components are independent stochastic
processes. Replacement this will always be determined by a probability measure $\mu$ on
$M$, equivalent to the volume measure, whereas the thinning follows a probability rule
on the Voronoi tessellation. We tacitly assume that all subsets of $M$ or $\mathscr{X}$
appearing in this paper are measurable.

Let a selection  function $S:\mathscr{C}\rightarrow \mathbb{R}^{+}$ be given, defined
on the space $\mathscr{C}$ of all possible Voronoi cells in $M$. If
$\vex{x}=(x_{1},\dots ,x_{k})$ is the current configuration, then in the next step
exactly one co-ordinate point $x_{J}$, $1\leq J \leq k$, is chosen at random and either
removed (in a thinning process) or reassigned (in a thinning-replacement process) to a
random position. The probability of choosing coordinate $x_{j}$ is proportional to the
value of the selection function $S(C_{x_{j}}^{\vex{x}})$ on the respective Voronoi cell
$C^{\vex{x}}_{x_j}$, i.e.
\begin{equation}\label{eq2}\Pr(J=j|\vex{x}):=\frac{S(C_{x_{j}}^{\vex{x}})}{\sum_{i=1}^{k}S(C_{x_{i}}^{\vex{x}})}.\end{equation}

There are many  quantifiable properties of a Voronoi cell upon which the function $S$
can be based. Some of these such as volume, perimeter or surface area or suitable high
dimensional generalisation, number of edges, faces, minimal or maximal internal angles
etc., are properties of any closed simplicial complex in the same underlying space.
Others, such as volume of the associated Voronoi flower (Voronoi flower of
$C_{x_{j}}^{\vex{x}}$ is the closure of the set $\left\{z\in M: C_{x_{j}}^{\vex{x}\cup
\left\{z\right\}}\neq C_{x_{j}}^{\vex{x}}\right\}$) or the distance to nearest
neighbour, are specific to the Voronoi structure but can be computed from the
properties of the individual cell. A further generalisation could allow $S$ to be a
function of $x_j$ and of $\vex{x}$, for example the weighted sum of cell functions in
the original sense with weights determined by the number of edges in the smallest
nearest-neighbour arc to $x_j$. In this sense, the degree of a selection function $S$
could be defined as the depth of the Voronoi nearest neighbour relation required to fix
its value. We will consider only the simplest kinds of first-degree selection
functions.  More specifically, we consider the following two classes of such processes.

(A) The ``volume-based",  or $v$-process. We require enough additional structure for
$(M,d)$ so that cell-volume can be defined. Let $M$ be a compact piecewise-smooth
manifold equipped with a measure $\lambda$, which is equivalent to the volume measure
on $M$ with a density bounded away from zero and from infinity. We assume that the
value of $S(C_{x_{j}}^{\vex{x}})$ is determined by the value of
$\lambda(C_{x_{j}}^{\vex{x}})$:
$$
S(C_{x_{j}}^{\vex{x}})=S_{v}(\lambda (C_{x_{j}}^{\vex{x}})).
$$
Without ambiguity we drop the subscript  in $S_{v}$ and consider
$S:\mathbb{R}^{+}\rightarrow\mathbb{R}^{+}$ to be a function of $\lambda
(C_{x_{j}}^{\vex{x}})$. If $S$ is increasing, then points with Voronoi cells of large
volume are more likely to be chosen to be culled or moved, and so the selection
pressure favours points with small cells, that is, points restricted by close
neighbours. A decreasing $S$ favours points with large cells, that is, isolated points
or points whose near neighbours fall within a limited arc. Functions of the form
$S(v)=v^{\alpha}$, $\alpha \in \mathbb{R}$ produce scale-independent models.

(B) The neighbour-based, or $n$-process.  The Voronoi tessellation determines for each
generator a set of its Voronoi-nearest neighbours. Formally, this set is defined as
$$
[x_{j}]^{\vex{x}}:=\left\{x_{i}:\vex{x}=(x_{1}\tdots x_{k})\, ,i \neq j,\, {\rm
card}(C_{x_{i}}^{\vex{x}}\cap C_{x_{j}}^{\vex{x}})>1 \right\}.
$$
We now assume that $S(C_{x_{j}}^{\vex{x}})=S_{n}({\rm card}([x_{j}]^{\vex{x}}))$,  for
some $S_n:\left\{1\tdots N\right\}\rightarrow \mathbb{R}$, where again we drop the
subscript where the context is clear. The $n$-process requires less structure on
$(M,d)$, but the selection function $S$ determines which, if any, types of cells are
favoured by the evolution.

Figures \ref{3oneD} to \ref{p6a5nprocs} below  illustrate some of the behaviours that
were observed in simulations.
\begin{figure}[h]
    \centering
        \includegraphics{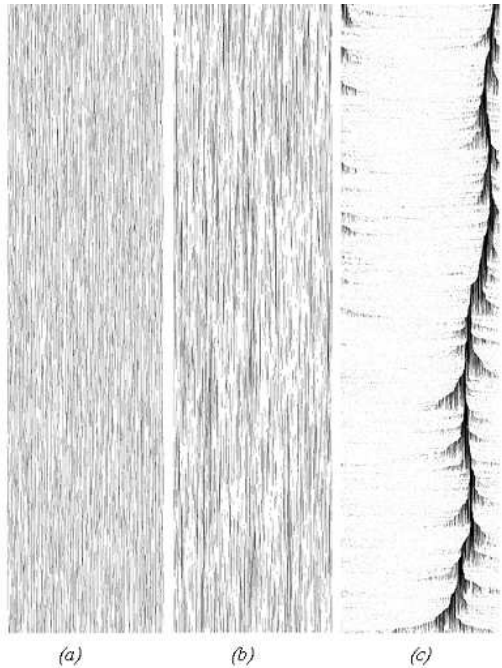}
        \caption{\footnotesize Evolution of $v$-processes on circle. Number of points, $N = 128$, number of steps, $T = 4096$, $S(v)=v^{\alpha}$, with (a) $\alpha = -1.0$, (b)  $\alpha = 0.5$, (c)  $\alpha=1.5.$}
    \label{3oneD}
\end{figure}

Fig.\,\ref{3oneD} depicts side-by-side realisations of three different $v$-processes on a circle with $S(v)=v^{\alpha}$, each having the same total number of points ($N=128$), but different $\alpha$ values. The base of each rectangle represents the circle $M$ opened out into a line segment by a cut, whereas the $y$-axis represents the time. The well defined clustering observed in $(c)$ was found to occur when $\alpha>1$ and $N$ is sufficiently large. The phase change was observable even for values of $\alpha$ close to 1, as seen in Fig.\,\ref{3phase1D}, where the time scale has been compressed by a factor of 5 from Fig.\,\ref{3oneD}, to bring it out more clearly.
\begin{figure}[h]
    \centering
        \includegraphics{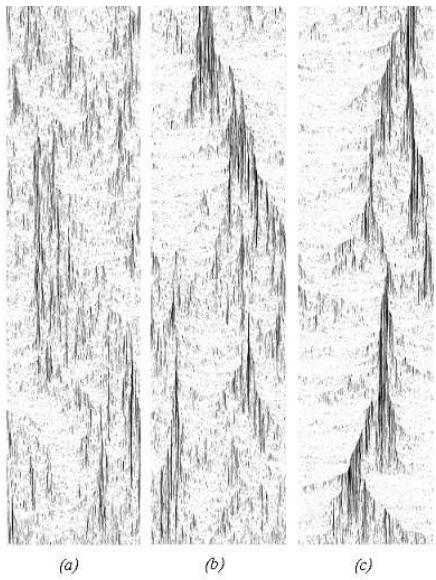}
        \caption{\footnotesize Phase change in the $v$-process on circle. $N = 128$, $T = 20480$, $S(v)=v^{\alpha}$, (a) $\alpha = 0.95$, (b)  $\alpha = 1.0$, (c)  $\alpha=1.05.$ }
    \label{3phase1D}
\end{figure}

Fig.\,\ref{vprocmatrix} shows the results of running four instances of the $v$-process on the unit square with $S(v)=v^{\alpha}$ and values of $\alpha$ ranging from -3.0 to 1.4. The phase change at $\alpha=1$ was observed just as in the one-dimensional case. For values of $\alpha \le 1$ a smooth gradation of degrees of clustering was observed without apparent interference from edge effects.

\begin{figure}[!h]
    \centering
        \includegraphics{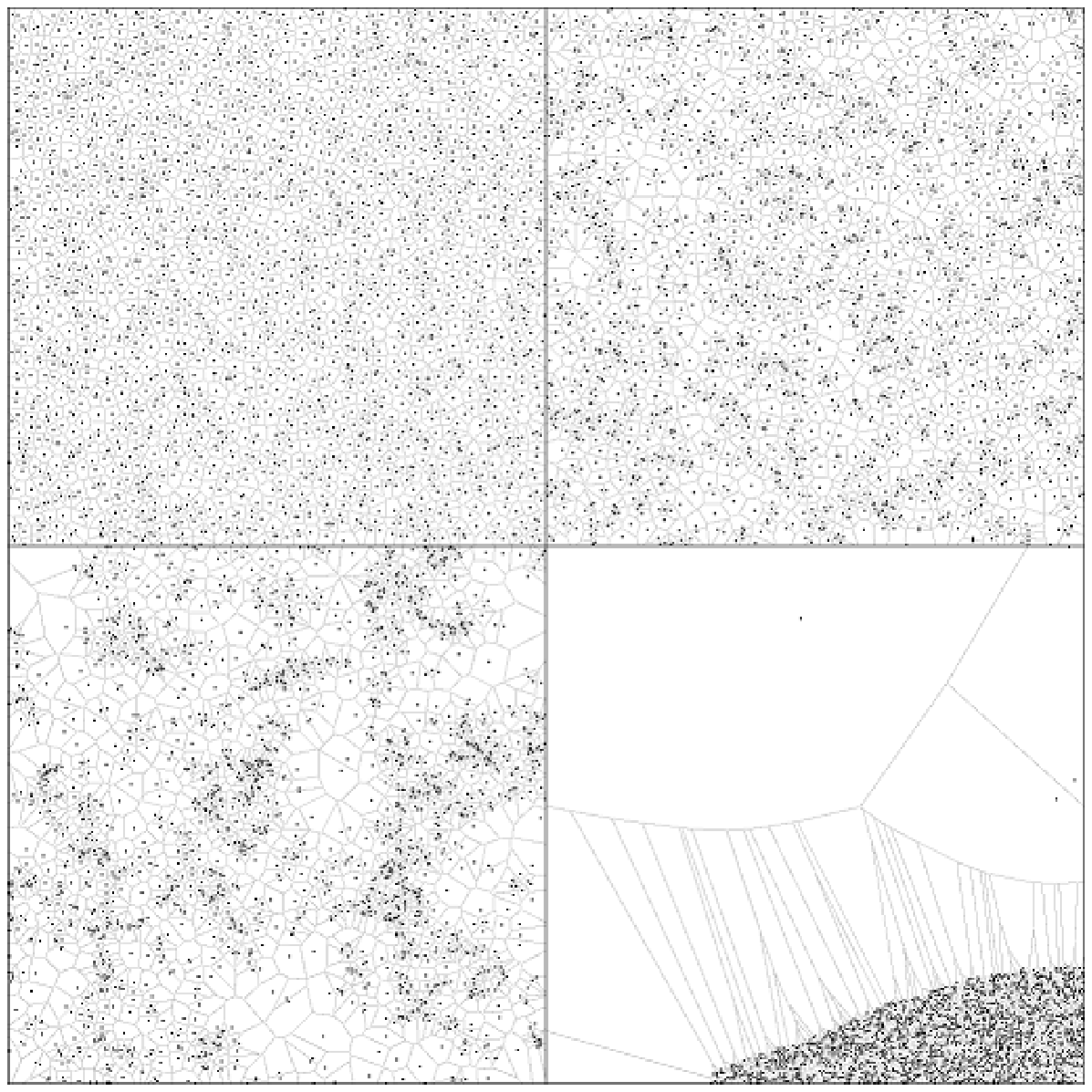}
        \caption{\footnotesize Realisations of the $v$-process with $S(v)=v^{\alpha}$ for $-3.0\leq  \alpha \leq 1.4$.
        Left-Right,Top-Bottom: $\alpha=$-3.0,  0.2, 1.0, 1.4. In each case there are $N=2000$ points with the number of cycles $T=24000$.}
    \label{vprocmatrix}
\end{figure}

Two instances of the $n$-process are shown in Fig.\,\ref{p6a5nprocs}.  The
$n$-processes produced a very rich collection of different clustering behaviours, but
in general it would seem that selections which favour cells with greater than average
numbers of neighbours (average is six) lead to less clustered configurations, and vice
versa.

\begin{figure}[!h]
    \centering
        \includegraphics{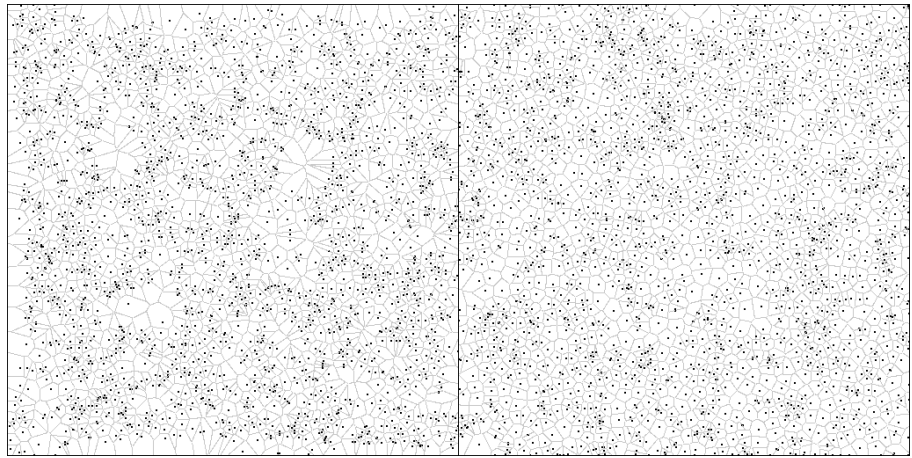}
        \caption{\footnotesize Realisations of the $n$-process with $N=2000$, $T=16000$. The pattern on the
        left was produced by the selection function $S(n)=0.1+|n-6|^2$ which ``favours" cells with six or close to six sides,
        that on the right by $S(n)=1,n\neq 5; S(5)=5000$ which strongly ``disfavours" cells with five sides. Note the apparent
        edge-effects in the former case. }
    \label{p6a5nprocs}
\end{figure}
\clearpage We will discuss the simulation results in more detail, and consider  the
problem of inferring the selection function from the statistics of a given point
pattern which has arisen from a $v$- or $n$-process elsewhere. In particular, we employ
Thiel's redundancy measure~\cite{OBS} and Baddeley and Van Lieshout's
$J$-function~\cite{BV96} as a first step in the characterisation of these patterns.

Simulations appear  to demonstrate ergodic-type behaviour for $v$-process with
$S(v)=v^{\alpha}$ when $\alpha \leq 1$ the behaviour has a very different character
which indicates that ergodicity may not take place. In the present paper we solve the
problem of ergodicity raised by these observations.

In Section 2 we prove,  under rather broad conditions, the ergodicity for both the $n$-
and $v$-processes. In Section 3 we present a crude explanatory model which explains the
phase change that is observed at $\alpha=1$ in the $v$-process with $S(v)=v^{\alpha}$.

\begin{rema}
Our process with replacement resembles a Gibbs sampler in that at each step only one
coordinate of the state is resampled. The main difference is that in a Gibbs sampler
the \emph{choice} of the coordinate to be resampled is deterministic, but the new
sample is from a distribution conditioned on the other current points. In our case, the
choice of point is from a distribution conditioned by the set of current points, and
the new sample is i.i.d. This latter choice could easily be modified so that the new
positions of the resampled coordinates are governed by a general set of conditional
distributions, as in the usual application of the Gibbs sampler.
\end{rema}

\section{Main results}

In this section we consider only processes with replacement. Assume $M$ is a compact
piecewise-$C^{2}$ manifold in $\mathbb{R}^{m}$, $m\ge 1$, with or without boundary,
endowed with the geodesic metric $d$. Let $\mu$ be a probabiloty and $\lambda$ a
measure on $M$, both equivalent to the volume measure $v$ on the manifold.

We define a discrete-time Markov $n$- or $v$-process with state  space
$\mathscr{X}=M^N$ as follows. Let $\vex{x} \in \mathscr{X}$, $A$ a Borel subset of
$\mathscr{X}$, $z\in M$, $\Delta\subseteq M$,
\begin{align*}
    \vex{x}[j,z]&:=(x_1\tdots x_{j-1},z,x_{j+1}\tdots x_N),
    \\A_{j}[\vex{x}]&:=\left\{z\in M:\vex{x}[j,z]\in A\right\} \text{ \ and }
    \\A[j,\Delta]&:=\left\{\vex{x}[j,z]\in \mathscr{X}: \vex{x}\in A,z\in \Delta \right\}.
\end{align*}

We also write, with a slight abuse of notation, $ x \in \vex{x}$  when
$\vex{x}=(x_1\tdots x_N)$ and $x=x_j$ for some $j, 1\leq j \leq N$.

The transition probability function is given by
\begin{equation}
    \label{eq3}P(\vex{x},A):=\Pr(\mathbf{X}_{n+1}\in A|\mathbf{X}_n=\vex{x})=\sum_{j=1}^{N}\mu
  (A_{j}[\vex{x}])\Pr(J=j|\vex{x}),
\end{equation}
where $\Pr(J=j|\vex{x})$ is the selection probability function \eqref{eq2}.

\begin{theo}
\emph{(i)} The $n$-process with replacement with a selection function $S:\left\{1,\dots
N\right\}\rightarrow \mathbb{R}^{+}$ which takes only positive values is Harris
ergodic.

\emph{(ii)} The $v$-process with replacement with a selection function
$S:(0,\lambda(M)]\rightarrow\mathbb{R}^{+}$, such that both $S$ and $1/S$ are locally
bounded, is Harris ergodic.
\end{theo}

\begin{rema}
It's not difficult to see that the assertion of part (i) holds  under more general
conditions as well.
\end{rema}

\begin{proof}
(i) Denote by $\mu^{\otimes N}$  the $N$-fold product measure on $\mathscr{X}$, and by
$P^{\left(k\right)}$, the $k$-step transition function generated by $P$. We start with
the following lemma.

\begin{lemo}
The n-process with replacement,  with $M$ and $\mu$ as specified above, satisfies the
following D$\ddot{o}$blin condition: there exist $\epsilon <1,\gamma >0$ such that if
$\muN(A)>\epsilon$ for a measurable $A\subseteq \mathscr{X}$, then
$P^{(N)}(\vex{x},A)\geq \gamma$  for any $\vex{x}\in \mathscr{X}.$
\end{lemo}

\begin{proof}
Let $B(z,\delta)=\left\{x\in M:d(x,z)<\delta\right\}$ be the $\delta$-ball around the
point $z\in M$, and let $H=\prod\limits_{i=1}^{N} B(z_i,\delta_i)\subset
\mathscr{X},\hspace{2mm} \text{ for } z_i\in M, \delta_i >0 $. It suffices to find a
constant $k_0 >0$, independent of $\left\{z_i,\delta_i: 1\leq i \leq N\right\}$, such
that for any $\vex{x}\in \mathscr{X}$
\begin{equation}
    \label{eq4}P^{(N)}(\vex{x},H) > k_0\muN(H).
\end{equation}
Indeed, since any open subset of $\mathscr{X}$ can be approximated in $\mu^{\otimes
N}$-measure arbitrarily well by disjoint unions of basic sets $H$, then by the
regularity of the measures on $\mathscr{X}$ the relation \eqref{eq4} implies that
\begin{equation}
    \label{eq5}
    P^{(N)}(\vex{x},A) > k_0\muN(A)
\end{equation}
for any $\vex{x}\in \mathscr{X}$ and $A \subseteq \mathscr{X}$.

First, note that for any $\vex{x}=(x_1,x_2\tdots x_N)\in \mathscr{X}$ and $\Delta
\subseteq M$  by \eqref{eq2} and  \eqref{eq3}, \begin{equation}
    \label{eq5.1}P\bigl(\vex{x},\left\{\vex{x}\right\}[j,\Delta]  \bigr) \geq k\mu(\Delta),
\end{equation}
 where
\begin{equation}
    \label{kbnd}
    k=\frac{\min_{j\leq N} S(j)}{N \max_{j\leq N} S(j)}>0.
\end{equation}

Therefore, for any $y\in M \text{ and } i\neq j$ we also have
$$
P(\vex{x}[i,y],\left\{\vex{x}[i,y]\right\}[j,\Delta])\geq k\mu(\Delta).
$$
Taking into account that the Markov chain can proceed from  $\vex{x}$ to the set
$\left\{\vex{x}\right\}[i,B(z_i,\delta_i)][j,B(z_j,\delta_j)]$ by first changing the
value of the $i^{\emph{th}}$ component and then that of the $j^{\emph{th}}$ one or vice
versa, we obtain from \eqref{eq5.1} that
\begin{multline*}
        P^{(2)}\bigl(\vex{x},\left\{\vex{x}\right\}[i,B(z_i,\delta_i)][j,B(z_j,\delta_j)]\bigr)
        \\=\int\limits_{B(z_i,\delta_i)}\int\limits_{B(z_j,\delta_j)}P\bigl(\vex{x},\left\{\vex{x}\right\}[j,dy_j]\bigr)P\bigl(\vex{x}
        [j,y_j],\left\{\vex{x}\right\}[j,y_j][i,dy_i]\bigr)
        \\+\int\limits_{B(z_j,\delta_j)}\int\limits_{B(z_i,\delta_i)}P\bigl(\vex{x},\left\{\vex{x}\right\}[i,dy_i]\bigr)P\bigl(\vex{x}
        [i,y_i],\left\{\vex{x}\right\}[i,y_i][j,dy_j]\bigr)
        \\\geq 2k^2\int\limits_{ B(z_i,\delta_i)}\mu(dy_i)\int\limits_{ B(z_j,\delta_j)}\mu(dy_j)=2k^2\mu(B(z_i,\delta_i))\mu(B(z_j,\delta_j)).
\end{multline*}
 Continuing in this way, we will clearly get
$$
P^{(N)}\bigl(\vex{x},\prod_{j=1}^{N}B(z_j,\delta_j)\bigr)\geq N!k^N\prod_{j=1}^{N}\mu (B(z_j,\delta_j)),
$$
which proves \eqref{eq4} and the lemma.

\end{proof}

The relation \eqref{eq5} clearly implies the irreducibilty of the \emph{n}-process.
Indeed,  \eqref{eq5} also implies that for any $\vex{x}\in \mathscr{X}$, $A\subseteq
\mathscr{X}$,
$$
P^{(N+1)}(\vex{x},A) > k\muN(A),
$$
and hence the process is aperiodic. Now the desired assertion is immediate from a
theorem in \cite{MT}, p.391. Part (i) of the theorem is proved.

(ii) Ergodicity of the $v$-process requires more analysis.  We wish to produce an
argument similar to that of Lemma 1. The measure $\mu$ determines probabilities for
points entering the configuration (the replacement points) and it will be sufficient
for our arguments to note that for any fixed $\delta$,
$$
\mu_{\delta}:= \inf\left\{\mu(B(x,\delta)): x \in M \right\}>0
$$
as a consequence of the smoothness and compactness of $M$. Similarly we can define,
$$
\lambda_{\delta}:= \inf\left\{\lambda(B(x,\delta)): x \in M \right\}>0.
$$

For the $v$-process we do not have a simple analogy to the lower bound derived from
\eqref{kbnd}. Here what is required is a positive lower bound for
$\lambda(B(x,\delta)\cap C^{\vex{x}}_{x})$,  which does not exist if $\vex{x}$ is
unrestricted. For example, if $x\in \partial M$, then a single other point in the
configuration $\vex{x}$ can cause the Voronoi cell around $x$ to be arbitrarily small
in volume and, consequently, the value of the selection function may become arbirarily
small on this cell. So first we need to define those points which are not too close to
the boundary.

Let $M_{\delta}:=\left\{x\in M: d(x,y)>\delta \text{ for any }y\in
\partial M\right\}$ for a $\delta>0$; if $M$ is without a boundary, then
$M_{\delta}:=M$. The smoothness of $M$ guarantees that $\mu(M_{\delta})>0$ and
$\lambda(M_{\delta})>0$ for $\delta$ small enough. Take any $x\in M_{\delta}$  and
$\vex{x}\in \mathscr{X}$, and define
$$
n(x,\vex{x},\delta):=\left\{\begin{array}{ll}{\rm card}(\left\{y\in
\vex{x}:0<d(x,y)<2\delta\right\})\quad &\text{for }x\in \vex{x};\\ N
&\text{otherwise}.\end{array}\right.
$$
So if $\vex{x}$ is such that $n(x,\vex{x},\delta)\leq 1$, then $B(x,\delta)$ intersects
with at  most two of the Voronoi cells $C^{\vex{x}}_{x_j}$. This means that if $x\in
M_{\delta}$ is fixed, then for all $\delta>0$ sufficiently small, and any $\vex{x}\in
\mathscr{X}$ such that $n(x,\vex{x},\delta)\leq 1$, we have, via the smoothness of $M$,
\begin{equation}
    \label{lamdel}
    \lambda (B(x,\delta)\cap C^{\vex{x}}_{x})\geq \lambda_{\delta/3}.
\end{equation}

Define for $\delta >0$,
$$
\mathscr{D}_{\delta}:=\left\{(x_1\tdots x_N)\in \mathscr{X}: \text{all }x_i\in M_{\delta}\text{ and }d(x_i,x_j)>2\delta,  i \neq j\right\}.
$$
Clearly we can choose $\delta$ small enough so that $\mu^{\otimes N}(\mathscr{D}_{\delta})>0$.

Now define a new measure $\phi$ on $\mathscr{X}$ as follows:
$$
\phi(A):=\mu^{\otimes N}(A\cap \mathscr{D}_{\delta}).
$$

\begin{lemo}
\emph{(i)} There are $\delta>0$ and $\epsilon>0$ such that for  every $\vex{x}\in
\mathscr{X}$,
\begin{equation}
\label{eq7}
P^{(N)}(\vex{x},\mathscr{D}_{\delta})>\epsilon.
\end{equation}

\emph{(ii)} There is a $\gamma>0$ such that for $\delta$ as  above and any $\vex{x}\in
\mathscr{D}_{\delta}$, and any Borel subset $A\subseteq \mathscr{X }$
\begin{equation}\label{eq8}P^{(N)}(\vex{x},A)\geq \gamma\phi(A).\end{equation}
\end{lemo}

\begin{proof}
(i) Let $\delta_0>0$  be small enough that $M_{\delta_0}$ has postive $\mu$ (and hence
positive $\lambda$) measure. By appealing to the smoothness structure of $M$ or
otherwise we can find $N+1$ disjoint open subsets of $M_{\delta_0}$. Denote these by
$L_1,L_2\tdots L_{N+1}$ and fix them for the rest of the proof. Let
$$
\delta_i:=\sup{\left\{\delta: \text{there are } N \text{ disjoint balls } B(z_k, 3\delta)\subset L_i\right\}},\quad i=1\tdots N+1.
$$
Clearly,
$$
\delta := \min_{0\leq i\leq N+1}\delta_i>0.
$$
Now for any $\vex{x}=(x_1\tdots x_N)\in \mathscr{X}$, at  least one of the $L_i$ will
be disjoint from $\left\{x_1\tdots x_N\right\}$. Therefore for such an $i$ we can find
$z^{(i)}_1\tdots z^{(i)}_N \in L_i$ such that

$$
B(z^{(i)}_j,3\delta)\cap B(z^{(i)}_k,3\delta)=\varnothing \quad \text{ for }j\neq k,
$$
$$
x_k \notin B(z^{(i)}_j,3\delta) \text{ for all }k,j \in \left\{1\tdots N\right\}.
$$
Obviously, $\prod_{j=1}^{N}B(z^{(i)}_j,\delta)\subset \mathscr{D}_{\delta}.$

As in Lemma 1, we consider a possible  sequence of $N$ steps of the Markov process in
each of which exactly one of the co-ordinate points $x_m$ of $\vex{x}$ moves into a
unique $B(z^{(i)}_{j'},\delta), j=1\tdots N$. At each step, the probability of the
chosen point landing in a $B(z^{(i)}_{j'},\delta)$ is greater than or equal to
$\mu_\delta$, so we need only be concerned with the probability of all $N$ original
co-ordinate points moving in $N$ consecutive steps of the process. This will have a
lower bound independent of the initial state $\vex{x}$ if we can find a bound $\kappa
>0$ with the following property: whenever we are at a stage at which some, but not all,
of the $x_m$ have moved and are in corresponding sets $B(z^{(i)}_{j'},\delta)$, there
is a $j$ such that $x_j$ has not moved yet, and
$$
\Pr(J=j|\vex{x}')>\kappa\text{, where } \vex{x}'\text{ is the current state of the process.}
$$

Let us assume that in the first $k\geq 1$  steps $k$ different original points
$x_{i_1}\tdots x_{i_{k}}$ have moved. We denote their new positions by $x'_{i_1}\tdots
x'_{i_{k}}$. Write $\nu :=\left\{i_1\tdots i_{k}\right\}$, $\xi:=\left\{1\tdots
N\right\}\backslash \nu = \left\{j: x_j \text{ has not moved} \right\}$ and
$$
S_{\Sigma}:=\sum_{j=1}^{N}S(\lambda(C^{\vex{x}'}_{x_j}))=\sum_{j\in \nu}+\sum_{j \in \xi}=:S_{\Sigma}^{\nu}+S_{\Sigma}^{\xi}.
$$
Then $p:=S_{\Sigma}^{\xi}/S_{\Sigma}= \Pr(J \in \xi |\vex{x}')$ is  the probability
that at the next step one of the remaining points will move.

If $p>1/2$, then for some $j\in \xi$,
\begin{equation}
\label{eq9}
\Pr(J=j|\vex{x}')\geq \frac{1}{2{\rm card}(\xi)} \geq \frac{1}{2N},
\end{equation}
otherwise
\begin{equation}
\label{eq10}
S_{\Sigma}^{\xi}\leq S_{\Sigma}^{\nu}.
\end{equation}

We only have to consider the latter case. For an $x'_j$, $j\in \nu$, we have
$$
x'_j \in B(z^{(i)}_{j'},\delta)
$$
for some $j'$, and so not only do we have
$$
d(x'_j, x'_m)>2\delta,
$$
when $j\neq m\in \nu$, but also for $m\in \xi$,  since $B(z^{(i)}_{j'},3\delta)$ was
free of any $x_k$ from the original state $\vex{x}$. This implies that for the set $W:=
\bigcup_{j\in \xi}B(x_j,\delta)$ one has $C^{\vex{x}'}_{x'_i}\cap W =\varnothing$ for
$i\in \nu$, and therefore $$ W\subset \bigcup_{j\in \xi}C^{\vex{x}'}_{x_j}.
$$
Hence,
$$
\max_{j\in \xi} \lambda(C^{\vex{x}'}_{x_j})\geq  \frac{\lambda(W)}{{\rm card}(\xi)}\geq
\frac{\lambda_{\delta}}{{\rm card}(\xi)}\geq \frac{\lambda_{\delta}}{N},
$$
while for all $j \in \nu$,
$$
\lambda (C^{\vex{x}'}_{x'_j})\geq \lambda_\delta.
$$
Due to the local boundedness of  $S$ and $S^{-1}$,
$$
0<b:=\inf_{\lambda_{\delta}/N\leq v \leq \lambda(M)}S(v)\leq \sup_{\lambda_{\delta}/N\leq v \leq \lambda(M)}S(v)=:B<\infty.
$$

Therefore by \eqref{eq2} and \eqref{eq10}
\begin{equation}
    \label{eq11}
    \max_{j \in \xi}\Pr(J=j|\mathbf{\vex{x}'}) \geq \frac{b}{S_{\Sigma}^{\xi}+{\rm card}(\nu)B}\geq \frac{b}{2NB}=:\kappa>0
\end{equation}
Clearly $\kappa \leq \frac{1}{2N}$, so from \eqref{eq9} and \eqref{eq11} we have,  by a
similar argument to that in the proof Lemma 1, that
\begin{equation}
    \label{eq112}
    P^{(N)}(\vex{x},\mathscr{D}_{\delta})\geq N!\kappa^{N-1} \mu_\delta^N,
\end{equation}
which proves \eqref{eq7}.

(ii) Let $A\subset \mathscr{X}$ and assume that $\phi(A)>0$. As in the proof of Lemma
1,  it will suffice to restrict attention to basic open sets
$$
H:=\prod_{i=1}^{N} B(w_i,r_i)\subset \mathscr{D}_{\delta},\qquad w_i\in M_{\delta},
$$
and as in the proof for (i), we specify a possible sequence of $N$  moves in the Markov
process and aim to produce a bound analogous to \eqref{eq112} for $P^{(N)}(\vex{x},H)$.
The term $\mu_\delta^N$ in \eqref{eq112} is easily seen to be replaceable in that case
by
$$
\prod_{i=1}^{N} \mu(B(w_i,r_i)=\phi(H),
$$
since it represents the product of the probabilities  that the points move \emph{into}
the target sets; the value $\kappa$ from (i) is no longer valid however, because,
retaining the notation above, we don't necessarily have that $d(x_j,x'_m)>2\delta$ when
$j\in \xi$ and $m \in \nu$. Thus, if we can find a new constant $\kappa' >0$ such that
at every intermediate state $\vex{x}'$
\begin{equation}
    \label{eq13}
    \max_{j \in \xi}{\Pr(J=j|\vex{x}')}> \kappa',
\end{equation}
we will have
\begin{equation}
    \label{eq14}
    P^{(N)}(\vex{x},H)>N!\kappa'^{N-1}\phi(H),
\end{equation} and \eqref{eq8} will have been proven.

As before, let $\vex{x}'$ be the state of the  process after $k$ steps, and let $j\in
\xi$. Since $\vex{x} \in \mathscr{D}_{\delta}$, for any other $i\in \xi$, we have
$d(x_j,x_i)>2\delta$. On the other hand, for all $l,m \in \nu$ we also have
$d(x'_l,x'_m)>2\delta$ as these points will be part of the configuration that we will
get after $N$ steps, which has to belong to $H \subset \mathscr{D}_{\delta}$. Let $m\in
\nu$ be such that
$$
d(x'_m,x_j)=\min_{i\in \nu}{d(x'_i,x_j)}.
$$
If $d(x'_m,x_j)<\delta/2$ then for  all other $i\in \nu$, $d(x'_i,x_j)>3\delta/2$, and
so the ball $B(x_j,3\delta/4)$ is contained entirely within $C^{\vex{x}'}_{x_j}\cup
C^{\vex{x}'}_{x'_m}$, and hence for sufficiently small $\delta>0$, by \eqref{lamdel},
$$
\lambda(C^{\vex{x}'}_{x_j})\geq \lambda_{\delta/8}.
$$
Alternately, $d(x'_i,x_j)\geq \delta/2$ for all $i \in \nu$, and so
$B(x_j,\delta/4)\subseteq C^{\vex{x}'}_{x_j}$.  In any case, $
\lambda(C^{\vex{x}'}_{x_j})\geq \lambda_{\delta/8}, $ and hence,
$$
\Pr(J=j|\vex{x}')\geq \frac{\lambda_{\delta/8}}{S_{\Sigma}}.
$$
We also have,
\begin{align*}
    S_\Sigma &\leq {\rm card}(\xi)\sup_{v\in [\lambda_{\delta/8},\lambda(M)]} S(v)+{\rm card}(\nu)\sup_{v\in [\lambda_\delta,\lambda(M)]}S(v)\\
    &\leq N\sup_{v\in [\lambda_{\delta/8},\lambda(M)]} S(v)=:Q <\infty.
\end{align*}
So we can take  $\kappa':=Q^{-1}\inf_{v\in [\lambda_{\delta/8},\lambda(M)]}S(v)$ in
\eqref{eq13}. Lemma 2 is proved.
\end{proof}

Following \S\hspace{5pt}5.4.3 of \cite{MT}, aperiodicity for the $v$-process is defined
as follows. Let
\begin{displaymath}
\begin{split}
 E_{\mathscr{D}_{\delta}}:=\left\{n\geq1: \text{there is a }\gamma_{n}>0 \text{ s.t. }
 P^{(n)}(\vex{x},A)\geq \gamma_{n}\phi(A)\right.\\\left.\text{ for any }\vex{x}\in \mathscr{D}_{\delta},A\subset \mathscr{X}\right\}.
\end{split}
\end{displaymath}
The $v$-process is called aperiodic if $\text{g.c.d.}( E_{\mathscr{D}_{\delta}})=1$.

\begin{lemo}
The $v$-process is aperiodic.
\end{lemo}

\begin{proof}
It is enough to find a $\gamma_{N+1}>0$ such that for any $\vex{x}\in
\mathscr{D}_{\delta}$, and any Borel subset $A\subseteq \mathscr{X}$
\begin{equation}
\label{eq12}
P^{(N+1)}(\vex{x},A)\geq \gamma_{N+1}\phi(A).
\end{equation}
From Lemma 2(ii)
\begin{multline*}
P^{(N+1)}(\vex{x},A)=\int_{\mathscr{X}}P(\vex{x},d\vex{y})P^{(N)}(\vex{y},A)\\\geq \int_{\mathscr{D}_{\delta}}P(\vex{x},d\vex{y})P^{(N)}(\vex{y},A)
\geq P(\vex{x},\mathscr{D}_{\delta})\gamma \phi(A).
\end{multline*}
Without loss of
generality we can choose $\delta$ so small that for any $N$ points $z_i\in M_{\delta},
1\leq i \leq N$, such that the balls $B(z_i,\delta)$ are disjoint, we can find a
$z_{N+1} \in M_{\delta}$ with
$$
B(z_{N+1},\delta)\cap \bigcup_{1\leq i \leq N}B(z_i,\delta)=\varnothing.
$$
For such a $\delta$ it is clear from \eqref{eq3} and the definition of $\mu_{\delta}$
that $P(\vex{x},\mathscr{D}_{\delta})>\mu_{\delta}$, which means that we can take
$\gamma_{N+1}:=\mu_{\delta}\gamma >0.$
\end{proof}

Now the assertion of part (ii) of the theorem is an immediate consequence of Lemmas 2
and 3.  Theorem 1 is proved.
\end{proof}

\section{Local behaviour --- a crude explanatory model}
\label{simplistic}

Simulations of the $v$-process with $S(v)=v^{\alpha}$ show avalanche-scale clustering
leading to the formation of a permanent tight cluster when $\alpha >1$ and $N$ is
sufficiently large, and a weak variable clustering when $\alpha \leq 1$. It would be of
interest to obtain some insight into the causes of this phase-change type of
phenomenon. In this subsection we will present a simplistic model at a physicist's
level of rigour which explains why such a transition occurs at the threshold value
$\alpha =1$.

Let $A \subset M$ be some small connected  ``test region" and $N_A(\vex{x}):={\rm
card}(\left\{i:x_i\in A\right\}).$ We consider the evolution of $N_A(\mathbf{X}_n)$.
Let
$$
\mathscr{S}_{B}(\vex{x}):=\sum_{x_i\in B}S(\lambda(C^{\vex{x}}_{x_i})), B\subset M.
$$
Given $\mathbf{X}_n = \vex{x}$, at the following step of the process the probability of
a point being lost from $A$ is
\begin{equation}
\label{eq17a} \frac{\mathscr{S}_A(\vex{x})} {\mathscr{S}_M(\vex{x})},
\end{equation}
while the probability of a fresh point entering $A$ is $\mu(A)$. As the total number
$N$ of points in the process  is typically much larger than $N_A(\mathbf{X}_n)$, it is
natural to expect that, due to the effect of a law of large numbers, the relative
fluctuations in $\mathscr{S}_M(\mathbf{X}_n)$ will be relatively small compared to
those in $\mathscr{S}_A(\mathbf{X}_n)$ (this is borne out by the results of
simulations). So let's assume, for simplicity, that
$\mathscr{S}=\mathscr{S}_M(\mathbf{X}_n))$ is constant.

We wish to find an approximation to \eqref{eq17a} as a function of
$N_A=N_A(\mathbf{X}_n).$ Another simplifying assumption (also supported by simulations)
is that the conditional (given $N_A$) distribution of cell volumes for cells whose
generators lie in $A$ is the same, modulo scale, for different values of $N_A$. That
is, for the conditional (on $N_A$) distribution function of the volume $V$ of a
randomly chosen cell with generator in $A$ one has
$$
\Pr(V\leq v|N_A)=g(v/m_A), \quad v>0,
$$
for some $g$, where  $m_A=\mathbf{E}(V|N_A)$. Furthermore, we can take
\begin{equation}
    \label{eq17}m_A=N_A^{-1}\mathbf{E}\sum_{X_{n,i} \in A}\lambda(C^{\mathbf{X_n}}_{X_{n,i}}|N_A)\approx\beta \lambda(A)/N_A,
\end{equation}
where $\beta$ is a quantity dependent only on the geometry of $A$ and the order of
magnitude of $N_A$, which reflects the fact that the union of the cells with generators
in $A$ overlaps $A$ itself. A rough calculation shows that $\beta = 1+O(N_A^{-1})$.
When $A$ is a square or a circle and $N_A$ is not often less than 20, $\beta$ can be
considered to be in the range $(1.0,1.5]$. Appealing to the law of large numbers, we
could write
$$
\mathscr{S}_A(\mathbf{X}_n)\approx N_A \mathbf{E}\bigl(S(\lambda(C^{\mathbf{X}_n}_{X_{n,I}}))\bigl|N_A\bigr)
=N_A\int_{0}^{\lambda(M)}S(v)dg(v/m_A),
$$
where $I$ stands for the index of a ``typical" cell with a generator $X_{n,I} \in A$.
Now since $S(v)=v^{\alpha}$, using \eqref{eq17} the integral above becomes
\begin{align*}
    \int_{0}^{\lambda(M)}v^{\alpha}dg(v/m_A)&=m_A^{\alpha}\int_{0}^{\lambda(M)/m_A}w^{\alpha}dg(w)\\
    &\approx N_A^{-\alpha}(\beta\lambda(M))^{\alpha}\int_{0}^{\infty}w^{\alpha}dg(w),
\end{align*}
making the natural assumption  that the last integral converges. Combining these
approximations we get the following estimate for the probability that a point is
removed from $A$ in one step:
$$
\mathscr{S}_A(\mathbf{X}_n)\mathscr{S}^{-1} \approx
N_A^{1-\alpha}(\beta\lambda(M))^{\alpha}\mathscr{S}^{-1}\int_{0}^{\infty}w^{\alpha}dg(w)=KN_A^{1-\alpha},
$$
for some constant  $K$. Thus, if $\Delta N_A$ denotes the change in $N_A$ in one step
of the $v$-process, we have \begin{equation}\label{eq18}\mathbf{E}(\Delta
N_A|N_A)\approx\mu(A)-KN_A^{1-\alpha}.\end{equation}

As the right-hand side of  \eqref{eq18} is an increasing function of $N_A$ when
$\alpha>1$, in this case we have a positive feedback condition for the mean of the
number $N_A$ of points in our test region $A$. This means that the process is bound to
quickly leave the ``intermediate" range of states characterised by diffuse, roughly
uniform spatial point distributions --- an observation that is in agreement with the
simulation data. Note also that, once the ``destabilising mechanisms" have transformed
the point distribution to a single (or a few) tight cluster(s), the assumptions on
which the crude local model was based are no longer valid.

On the other hand, when $\alpha \leq 1$,  the relation \eqref{eq18} expresses either
neutral ($\alpha =1$) or negative feedback. Hence one expects ``local stability" from
the process behaviour: small clusters of points will form and disappear without any
``global" dramatic changes for the whole picture.

In such cases we can expect $N_A$  to take values close to $N\mu(A)$, so we can
estimate $K$ by the equation
$$
K(N\mu(A))^{1-\alpha}\approx \mu(A),
$$
from which we conclude that
$$
K\approx\mu(A)^{\alpha}N^{\alpha-1}, \quad \alpha \leq1.
$$
This approximation is reasonably well-supported by simulations.

\end{document}